\documentclass[11pt,a4paper]{article}
\usepackage[T1]{fontenc}
\usepackage[utf8]{inputenc}
\usepackage{lmodern}
\usepackage{amsmath,amssymb,amsthm,mathtools}
\usepackage[margin=29mm,headheight=14pt]{geometry}
\usepackage{microtype}
\usepackage{enumitem}
\usepackage{needspace}
\usepackage{xcolor}
\usepackage{hyperref}
\hypersetup{colorlinks=true,linkcolor=blue!45!black,citecolor=blue!45!black,
  urlcolor=blue!45!black,pdftitle={On Fujita's Advanced Generalized Mukai Conjecture}}
\setlist[enumerate,1]{label=\textup{(\roman*)},leftmargin=2.6em,itemsep=4pt,topsep=5pt}
\allowdisplaybreaks[1]
\newtheorem{theorem}{Theorem}[section]
\newtheorem{proposition}[theorem]{Proposition}
\newtheorem{lemma}[theorem]{Lemma}
\newtheorem{corollary}[theorem]{Corollary}
\theoremstyle{remark}
\newtheorem{remark}[theorem]{Remark}
\numberwithin{equation}{section}
\newcommand{\PP}{\mathbf P}
\newcommand{\QQ}{\mathbf Q}
\newcommand{\CC}{\mathbf C}
\newcommand{\RR}{\mathbf R}

\newcommand{\OO}{\mathcal O}
\newcommand{\Mbar}{\overline{\mathcal M}}
\newcommand{\pt}{\mathrm{pt}}
\newcommand{\vir}{\mathrm{vir}}
\newcommand{\NE}{\overline{\mathrm{NE}}}
\newcommand{\AGM}{\mathrm{AGM}}
\DeclareMathOperator{\Locus}{Locus}
\DeclareMathOperator{\RatCurves}{RatCurves}
\DeclareMathOperator{\ev}{ev}
\DeclareMathOperator{\pr}{pr}
\DeclareMathOperator{\Supp}{Supp}

\DeclareMathOperator{\End}{End}

\title{On Fujita's Advanced Generalized\protect\\ Mukai Conjecture}
\author{Duo Li}
\date{}

\begin{document}
\maketitle
\vspace{-1.7em}
\begin{abstract}
We prove Fujita's advanced generalized Mukai conjecture for smooth complex
Fano varieties. Our starting point is the existence of independent two-point
Gromov--Witten correspondences established by OpenAI. Under the inequality
$\rho_X(\iota_X-1)\geq\dim X-1$, their anticanonical degrees are less than
$2\iota_X$. We show that they determine finitely many unsplit families of
rational curves whose chains connect $X$. We use these families to bound
the anticanonical degrees throughout the successive minimal horizontal
construction, so every family in that construction is unsplit. Fujita's
contraction theorem and inductive classification then give the required
list of varieties.
\end{abstract}

\section{Introduction}

Throughout the paper, $X$ is a smooth connected complex projective Fano
variety of positive dimension $n$. We denote its Picard number by $\rho_X$
and its pseudoindex by
\[
 \iota_X=\min\{-K_X\cdot C:\ C\subset X
       \text{ is an irreducible rational curve}\}.
\]
The Fano index $j_X$ is the largest positive integer such that
$-K_X\sim j_XL$ for a line bundle $L$. Kobayashi and Ochiai
\cite{KO} proved that $j_X\leq n+1$, with equality precisely for
$\PP^n$. Mukai made a similar  conjecture \cite{Mukai} for varieties of higher 
Picard numbers:
\[
 \rho_X(j_X-1)\leq n,
 \qquad
 \rho_X(j_X-1)=n\ \Longleftrightarrow\
 X\simeq(\PP^{j_X-1})^{\rho_X}.
\]
Following Wi\'sniewski's work on the pseudoindex \cite{Wis},
Bonavero--Casagrande--Debarre--Druel \cite{BCDD} formulated the generalized
Mukai conjecture, in which $j_X$ is replaced by $\iota_X$. And they  proved it upto dimension four.

Andreatta--Chierici--Occhetta \cite{ACO} proved the five-dimensional case.
Novelli and Occhetta \cite{NO} established the conjecture for
$\iota_X\geq(n+3)/3$, and Novelli \cite{Novelli} extended this range to
$\iota_X>n/3$. The conjecture is also known for $\QQ$-factorial Gorenstein
toric Fano varieties by Casagrande \cite{Casagrande}, for
$\QQ$-factorial Gorenstein horospherical Fano varieties by Pasquier
\cite{Pasquier}, and for $\QQ$-factorial spherical Fano varieties by
Gagliardi--Hofscheier--Pearson \cite{GHP}. Pearson \cite{Pearson} proved
the index version for Fano varieties whose canonical ambient toric
variety is smooth.

We use the following result of OpenAI.

\begin{theorem}[{\cite[Theorem 1.1]{OpenAI}}]\label{thm:GM}
Every smooth complex projective Fano variety $X$ of positive dimension
satisfies
\[
 \rho_X(\iota_X-1)\leq\dim X.
\]
Equality holds if and only if $X\simeq(\PP^{\iota_X-1})^{\rho_X}$.
\end{theorem}

Fujita \cite[Conjecture 1.4]{Fujita} proposed a classification one step
beyond the equality case, under the assumption
$\rho_X(\iota_X-1)\geq n-1$. His inductive results
\cite[Proposition 1.8 and Remark 1.9]{Fujita}, together with the
characterization of quadrics by Dedieu and H\"oring \cite{DH}, give this
classification when $\rho_X\leq3$ or $n\leq5$.

The proof of Theorem~\ref{thm:GM} provides the additional input we need.
For $\iota_X\geq2$, OpenAI \cite[Propositions 4.1 and 5.3]{OpenAI}
constructed $\rho_X$  positive stable maps whose numerical classes 
$\gamma_1,\ldots,\gamma_{\rho_X}$ are independent. And by using the ideas form the Gromov-Witten theory, OpenAI constructed  a nonzero ordered product of
two-point correspondences (which implies any two points can be connected by  chains from $\gamma_1,\ldots,\gamma_{\rho_X}$) and showed 
\[
 \sum_{j=1}^{\rho_X}(-K_X\cdot\gamma_j-1)\leq n.
\]
Under Fujita's numerical hypothesis, their anticanonical degrees are at most
$\iota_X+1<2\iota_X$. We use the support of the correspondences to obtain
unsplit rational-curve families whose chains contain all these numerical
directions. Their classes span $N_1(X)_{\RR}$, which forces the associated
rational-chain quotient to be a point. We then compare these families
with the minimal families in Fujita's successive quotient construction.
This comparison shows that every minimal family is unsplit, precisely the
hypothesis needed in \cite[Corollary 6.4]{Fujita}. The classification
follows from Fujita's inductive step \cite[Proposition 5.3(b)]{Fujita}.

We use the Grothendieck convention for projectivization: $\PP_Y(E)$
parametrizes one-dimensional quotients of $E$. For $m\geq2$, we write
\[
 B_m=\PP_{\PP^m}\bigl(\OO^{\oplus(m-1)}\oplus\OO(1)\bigr),
 \qquad I_m=\PP_{\PP^m}(T_{\PP^m}).
\]
Thus $I_m$ is the incidence variety of points and hyperplanes in $\PP^m$.
We denote a smooth $m$-dimensional quadric by $Q^m$ and interpret an empty
product as a point.

\begin{theorem}\label{thm:main}
Suppose that $\rho_X(\iota_X-1)\geq n-1$. Then $\iota_X\geq2$ and
\[
 \delta:=n-\rho_X(\iota_X-1)\in\{0,1\}.
\]
The following classification holds.
\begin{enumerate}
\item If $\delta=0$, then $X\simeq(\PP^{\iota_X-1})^{\rho_X}$.
\item If $\delta=1$ and $\rho_X=1$, then $n=\iota_X\geq3$ and
      $X\simeq Q^n$.
\item If $\delta=1$ and $\rho_X\geq2$, then $X$ is isomorphic to one of
\[
\begin{gathered}
 Q^{\iota_X}\times(\PP^{\iota_X-1})^{\rho_X-1}
       \qquad(\iota_X\geq3),\\
 B_{\iota_X}\times(\PP^{\iota_X-1})^{\rho_X-2},\qquad
 I_{\iota_X}\times(\PP^{\iota_X-1})^{\rho_X-2},\\
 \PP^{\iota_X}\times(\PP^{\iota_X-1})^{\rho_X-1}.
\end{gathered}
\]
\end{enumerate}
\end{theorem}

The restriction in the quadric case is necessary because
$Q^2\simeq\PP^1\times\PP^1$ has Picard number two.
Theorem~\ref{thm:auxiliary} below gives Fujita's formulation with an
auxiliary lower bound on $\rho_X$.

A related classification was recently obtained by Watanabe
\cite[Theorem 1.2]{Watanabe} for klt Fano varieties of nef complexity
less than two. The nef complexity is defined using decompositions of
$-K_X$ into nef Cartier divisors. Our numerical hypothesis involves the
pseudoindex, and our proving strategy is different.

We recall the correspondence theorem in Section~\ref{sec:correspondences}.
Section~\ref{sec:curves} contains the passage to unsplit families, and
Section~\ref{sec:classification} completes the classification.

\section{Independent curve correspondences}\label{sec:correspondences}

We work with integral numerical curve classes. The real spaces
$N^1(X)_{\RR}$ and $N_1(X)_{\RR}$ are dual under the intersection pairing
and have the same dimension $\rho_X$. For an integral numerical class $\beta$, we
write $d_\beta=-K_X\cdot\beta$ and denote by
$\Mbar_{0,m}(X,\beta)$ the stack of genus-zero stable maps with $m$
marking points and numerical degree $\beta$.

For this notation, we take the disjoint union of the usual stable-map
stacks over the integral homology classes with numerical image $\beta$.
Only finitely many nonempty components occur: an ample multiple of
$-K_X$ fixes the projective degree, and the corresponding stable-map
space is proper and of finite type \cite[Sections 2--5]{FP}.
We write $\beta>0$ when this degree is represented by a nonconstant
genus-zero stable map. If its nonconstant components map with degrees
$e_\nu$ onto reduced rational curves $C_\nu$, then
\begin{equation}\label{eq:degree-sum}
 \beta=\sum_\nu e_\nu[C_\nu],
 \qquad d_\beta=\sum_\nu e_\nu(-K_X\cdot C_\nu)
       \geq\iota_X\sum_\nu e_\nu.
\end{equation}

\subsection{The two-point operators}

Evaluation at the markings gives proper morphisms
$\ev_i:\Mbar_{0,2}(X,\beta)\to X$. The virtual fundamental class has
complex dimension $n+d_\beta-1$ \cite{Behrend}. Following
\cite[Section 2]{OpenAI}, we set
\[
 H=\bigoplus_{k=0}^{n}H^k,\qquad H^k=H^{2k}(X,\CC),
 \qquad (u,v)=\int_X u\cup v,
\]
and put $H^k=0$ outside this range. We denote by $\pt\in H^n$ the point
class (which is actually represented by a top form and we call it "point class" just because its dual  in  the Poincar\'e pairing is reprensented by a multiple of a point), normalized by $\int_X\pt=1$.

The two-point invariant and its associated operator are defined by
\begin{equation}\label{eq:operator}
 (S_\beta u,v)=\langle u,v\rangle_\beta
 :=\int_{[\Mbar_{0,2}(X,\beta)]^{\vir}}
             \ev_1^*u\cup\ev_2^*v,
 \qquad u,v\in H.
\end{equation}
The nondegenerate Poincar\'e pairing makes $S_\beta\in\End_{\CC}(H)$
well defined. Equivalently, $S_\beta$ is induced by the algebraic
correspondence
\[
 K_\beta=(\ev_1,\ev_2)_*[\Mbar_{0,2}(X,\beta)]^{\vir}
       \in A_{n+d_\beta-1}(X\times X)_{\QQ}.
\]
We will use its standard grading and symmetry properties
\cite[Section 2]{OpenAI}:
\begin{equation}\label{eq:grading}
 S_\beta(H^k)\subseteq H^{k+1-d_\beta},
 \qquad S_\beta^*=S_\beta.
\end{equation}
Here the adjoint is taken for the bilinear Poincar\'e pairing. The grading
comes from the virtual dimension, and the symmetry from exchanging the
two markings. 
\begin{remark}
    $S_\beta$ can be understood as follows: if $S_\beta (\pt)$ is nonzero, then the first evaluation morphism $\ev_1$ is surjective as $K_\beta$ is a class supported on $(\ev_1,\ev_2) \Mbar_{0,2}(X,\beta)$.
\end{remark}
\subsection{The external input}

We use the following consequence of the spectral and correspondence
results in \cite{OpenAI}. We index the degrees in their order of
application.

\begin{theorem}[{\cite[Propositions 4.1 and 5.3]{OpenAI}}]
\label{thm:independent}
If $\iota_X\geq2$, there are positive integral numerical stable-map
degrees $\gamma_1,\ldots,\gamma_{\rho_X}$, linearly independent in
$N_1(X)_{\RR}$, such that
\begin{equation}\label{eq:independent}
 S_{\gamma_{\rho_X}}\cdots S_{\gamma_1}\ne0,
 \qquad \sum_{j=1}^{\rho_X}(d_{\gamma_j}-1)\leq n.
\end{equation} Here $d_{\gamma_j}$ is the intersection number of $-K_X$ with $\gamma_j$.
\end{theorem}

\section{From correspondences to unsplit families}\label{sec:curves}

A family of rational curves on $X$ means an irreducible component
$V\subset\RatCurves^n(X)$ of the normalized space of reduced irreducible
rational curves. We write $[V]$ for its numerical class and
$\Locus(V)$ for the union of its members. The family is \emph{unsplit}
if $V$ is proper, \emph{dominating} if $\Locus(V)$ is dense in $X$,
and \emph{covering} if $\Locus(V)=X$.
An unsplit family has a proper universal family, so its locus is closed;
in particular, an unsplit dominating family is covering. We use the
associated Chow families as in \cite[II.2.11]{Kollar} and
\cite[Definitions 2.1--2.3]{Fujita}.

\begin{lemma}\label{lem:short-degrees}
Suppose that $\rho_X(\iota_X-1)\geq n-1$. Then $\iota_X\geq2$, and
the degrees in Theorem~\ref{thm:independent} satisfy
\[
 \delta=n-\rho_X(\iota_X-1)\in\{0,1\},
 \qquad \iota_X\leq d_{\gamma_j}\leq\iota_X+1<2\iota_X.
\]
After reversing their order if necessary, we have
\begin{equation}\label{eq:point-product}
 S_{\gamma_{\rho_X}}\cdots S_{\gamma_1}(\pt)\ne0.
\end{equation}
\end{lemma}

\begin{proof}
A smooth Fano curve is $\PP^1$ and has pseudoindex two. For $n\geq2$,
the hypothesis excludes $\iota_X=1$. We may therefore apply
Theorem~\ref{thm:independent}. We put
\[
 M=\sum_j(d_{\gamma_j}-1),\qquad
 T=S_{\gamma_{\rho_X}}\cdots S_{\gamma_1}.
\]
Every positive stable-map degree has anticanonical degree at least $\iota_X$ by
\eqref{eq:degree-sum}, so
\[
 n-1\leq\rho_X(\iota_X-1)\leq M\leq n,
 \qquad 0\leq\sum_j(d_{\gamma_j}-\iota_X)\leq\delta\leq1.
\]
Each summand is a nonnegative integer, which gives the asserted degree
bounds.

By \eqref{eq:grading}, the nonzero operator $T$ lowers codimension by
$M$. If $M=n$, its only possible nonzero graded component is
$H^n\to H^0$. Since $H^n=\CC\pt$, we obtain $T(\pt)\ne0$.
If $M=n-1$, the only possible nonzero components are
$H^n\to H^1$ and $H^{n-1}\to H^0$.
When $T(\pt)=0$, there is an $a\in H^{n-1}$ with
$T(a)=c\,1$ for some $c\ne0$. Thus
\[
 (a,T^*(\pt))=(T(a),\pt)=c\ne0.
\]
The self-adjointness of each factor gives
$T^*=S_{\gamma_1}\cdots S_{\gamma_{\rho_X}}$. Reversing the order
therefore gives \eqref{eq:point-product}.
\end{proof}

\begin{lemma}[Short degrees]\label{lem:unsplit}
For a positive degree $\beta$ with $d_\beta<2\iota_X$, every genus-zero
stable map of degree $\beta$ has a unique nonconstant component, mapped
birationally onto its reduced image. There are finitely many
irreducible components of $\RatCurves^n(X)$ of numerical class
$\beta$, and all are unsplit.
\end{lemma}

\begin{proof}
The standard short-degree criterion follows here directly from
\eqref{eq:degree-sum}: the positive integer $\sum_\nu e_\nu$ is less
than two. Hence it equals one. Every contracted component maps to a
point of the unique reduced rational image, so all marked evaluations
lie on that image.

We use the properness of the stable-map space and its cycle morphism to
the Chow variety \cite[Sections 2--5]{FP}. In the unmarked space
$\Mbar_{0,0}(X,\beta)$, no contracted component can occur: a contracted
tree attached to the unique nonconstant component would have a leaf
with only one special point, contrary to stability. We denote by \(Y_\beta\subset\operatorname{Chow}(X)\) the reduced
image of the cycle map
\[
\operatorname{cyc}\colon
\overline{\mathcal M}_{0,0}(X,\beta)
\longrightarrow \operatorname{Chow}(X),
\qquad
[f\colon C\to X]\longmapsto f_*[C].
\]
The preceding paragraph shows that the closed points of \(Y_\beta\)
correspond exactly to the reduced irreducible rational curves of
numerical class \(\beta\). Since
\(\overline{\mathcal M}_{0,0}(X,\beta)\) is proper and of finite type,
its cycle image \(Y_\beta\) is closed and proper and has finitely many
irreducible components \(Y_{\beta,1},\ldots,Y_{\beta,t}\), each endowed
with its reduced structure.

Numerical class is locally constant on
\(\operatorname{RatCurves}^{n}(X)\), so the locus of curves of class
\(\beta\) is a union of full irreducible components. Since \(Y_\beta\)
contains every such curve, the Chow-space description of
\(\operatorname{RatCurves}^{n}(X)\) identifies these components with
the normalizations \(V_i:=Y_{\beta,i}^{\nu}\)
\textup{[11, II.2.11]}.
Each normalization morphism \(V_i\to Y_{\beta,i}\) is finite.
As \(Y_{\beta,i}\) is proper, \(V_i\) is proper as well.
Thus there are only finitely many rational-curve components of
numerical class \(\beta\), and each is unsplit.
\end{proof}

\subsection{Chains detected by the correspondences}

The incidence argument in \cite[Lemma 6.3]{OpenAI} treats curves of
minimum anticanonical degree. We use the same support argument for
degrees below $2\iota_X$, keeping track of the starting point of the
chain.

\begin{proposition}\label{prop:chains}
Suppose that $\beta_1,\ldots,\beta_k>0$ satisfy
$d_{\beta_j}<2\iota_X$ and
$S_{\beta_k}\cdots S_{\beta_1}(\pt)\ne0$.
There are unsplit families $V_1,\ldots,V_k$, with $[V_j]=\beta_j$,
such that every $x\in X$ starts a chain of member curves
$C_j\in V_j$ with
\[
 x\in C_1,\qquad C_j\cap C_{j+1}\ne\varnothing\quad(1\leq j<k).
\]
\end{proposition}

\begin{proof}
We write $\mathcal M_j=\Mbar_{0,2}(X,\beta_j)$ and form the proper
fiber product
\[
 W=\mathcal M_1\mathbin{\times_X}\cdots\mathbin{\times_X}\mathcal M_k,
 \qquad e:W\longrightarrow X\times X.
\]
The fiber products identify the second evaluation on one factor with
the first evaluation on the next; $e$ records the two remaining
endpoints. We set $P=S_{\beta_k}\cdots S_{\beta_1}$.

The usual composition of algebraic correspondences gives a cycle
supported on $e(W)$ that acts on cohomology as $P$.
More explicitly, refined pullback along the matching diagonals in
powers of the smooth variety $X$ gives a class
\[
 \Theta=\Delta^!\bigl([\mathcal M_1]^{\vir}\times\cdots
                         \times[\mathcal M_k]^{\vir}\bigr)
       \in A_*(W)_{\QQ}.
\]
Here $\Delta:X^{k-1}\hookrightarrow(X\times X)^{k-1}$ is the product
of the diagonals, pulled back by the evaluation pairs
$(\ev_2^{(j)},\ev_1^{(j+1)})$; for $k=1$, the refined pullback is
the identity. The correspondence $e_*\Theta$ has the required action. This is the
standard compatibility of refined intersection with proper pushforward
and composition \cite[Sections 6.2 and 16.1]{Fulton}; it does not require
the virtual classes to be effective.

The first endpoint image $A=\pr_1(e(W))$ is closed. If $x\notin A$,
the correspondence is disjoint from $\{x\}\times X$, and its refined
intersection with this subvariety is zero. Compatibility with cycle
classes gives $P([x])=0$. Since every closed point represents
$\pt\in H^n$, this contradicts the hypothesis. We conclude that
$A=X$. By Lemma~\ref{lem:unsplit}, the reduced images of these stable
maps are rational curves containing both evaluation points.

It remains to choose the same tuple of families for every starting
point. For each $j$, Lemma~\ref{lem:unsplit} gives finitely many
unsplit families $V_{j,1},\ldots,V_{j,t_j}$ of class $\beta_j$.
For a fixed tuple of indices $\mathbf i=(i_1,\ldots,i_k)$, we take the
incidence space of chains whose $j$th member belongs to $V_{j,i_j}$.
This space is proper: it is obtained from the proper universal cycle
supports by marking two points on each member and matching consecutive
endpoints. Its first endpoint image $A_{\mathbf i}\subset X$ is
therefore closed. The preceding paragraph gives
\[
 X=\bigcup_{\mathbf i}A_{\mathbf i}.
\]
This is a finite union, so the irreducibility of $X$ implies that some
$A_{\mathbf i}$ equals $X$. The corresponding families have the
required property. In particular, $V_1$ is covering.
\end{proof}

\subsection{Rational-chain quotients}

For finitely many Chow families of connected rational one-cycles, we
use the equivalence relation generated by chains of their supports.
The associated quotient is a proper surjective morphism with connected
fibers
\[
 q:U\longrightarrow Z,
\]
where $U\subset X$ is dense open and $Z$ is quasiprojective. Each fiber
is an equivalence class for chains in $X$
\cite[IV.4.16]{Kollar}, \cite[Theorem 3.2]{Fujita}.
Each fiber is proper over $\CC$ and closed in $X$. Every member of
one of the specified families that meets the fiber lies entirely in it.
We use this saturation property below.

\begin{lemma}\label{lem:rc-connected}
Suppose that $V_1,\ldots,V_h$ are unsplit families whose classes span
$N_1(X)_{\RR}$. If every point of a dense open subset of $X$ starts a
chain of their members containing at least one member of each family,
then these families rationally chain connect $X$.
\end{lemma}

\begin{proof}
We denote the dense open subset in the hypothesis by $G$ and take the
rational-chain quotient $q:U\to Z$, with $Z$ integral.
Suppose that $\dim Z>0$. A hyperplane section of a quasiprojective
embedding of $Z$ gives a nonzero effective Cartier divisor $B$ with
nonempty complement. The Weil divisor closure
$D=\overline{q^*B}$ in $X$ is effective and nonzero. It is Cartier
because $X$ is smooth, and
\[
 \Supp D\cap U=q^{-1}(\Supp B).
\]
We choose $x\in G\cap U\setminus\Supp D$. The fiber
$F=q^{-1}(q(x))$ lies in $U$ and is disjoint from $D$.
By saturation, the entire chain furnished by the hypothesis lies in
$F$. For a member $C_j\in V_j$ from this chain, we therefore have
\[
 D\cdot[V_j]=D\cdot C_j=0\qquad(1\leq j\leq h).
\]
The spanning assumption gives $D\equiv0$. This contradicts
$D\cdot A^{n-1}>0$ for an ample divisor $A$, since $D$ is a nonzero
effective divisor.

Thus $Z$ is a point. Properness of $q$ makes $U$ proper over $\CC$;
its inclusion in the separated variety $X$ consequently has closed
image. Since $U$ is dense, we have $U=X$. The quotient theorem now
gives rational-chain connectedness on all of $X$.
\end{proof}

For a proper surjective morphism $q:U\to Z$ as above, a curve is
\emph{horizontal} if it meets $U$ and its image under $q$ is not a
point. Following \cite[Definition 2.2]{Fujita}, we call a family
\emph{horizontal dominating} if its locus dominates $Z$ and none of
its members meeting $U$ is contracted by $q$.

\begin{lemma}\label{lem:horizontal}
Suppose that the unsplit families $V_1,\ldots,V_h$ rationally chain
connect $X$. If $q:U\to Z$ is proper and surjective, with $U\subset X$
dense open and $Z$ an integral quasiprojective variety of positive
dimension, then one of these families is horizontal dominating for $q$.
\end{lemma}

\begin{proof}
We first show that a fixed family cannot contain both a horizontal
member $C_h$ and a member $C_v$ meeting $U$ that is contracted by $q$.
The fiber over the image $z_v$ of $C_v\cap U$ is proper and hence
closed in $X$, so $C_v\subset q^{-1}(z_v)\subset U$.
We choose two distinct points of $q(C_h\cap U)$, the first different
from $z_v$, and a hyperplane section $B$ of $Z$ passing through the
first and avoiding the second and $z_v$. As in the preceding proof,
the closure $D=\overline{q^*B}$ is an effective Cartier divisor on $X$.
It is disjoint from $C_v$, meets $C_h$, and does not contain $C_h$.
Consequently
\[
 D\cdot C_v=0<D\cdot C_h,
\]
contrary to the numerical equivalence of members of the same family.

We denote by $J$ the set of indices of families that contain a
horizontal member. For every $z\in Z$, there is a chain from a point
of $q^{-1}(z)$ to a point outside this fiber. The first member that
leaves the fiber is horizontal: otherwise its intersection with $U$
would map constantly to $z$, and closedness of the fiber would force
the whole curve to remain in it. Hence
\[
 Z=\bigcup_{j\in J}q\bigl(\Locus(V_j)\cap U\bigr).
\]
Each set on the right is closed, since $V_j$ is unsplit and $q$ is
proper. A finite union of proper closed subsets cannot cover the
irreducible variety $Z$, so one of these images equals $Z$.
The first paragraph shows that every member of the corresponding
family meeting $U$ is horizontal. This is the required family.
\end{proof}

\subsection{The minimal sequence}

We now use the construction of Novelli--Occhetta
\cite[Construction 4.3]{NO}, in the form recalled by Fujita
\cite[Construction 3.5]{Fujita}. We choose a dominating family $W_1$
of least anticanonical degree. If its associated Chow family does not
rationally chain connect $X$, we choose a horizontal dominating
family $W_2$ of least degree for the quotient. We continue with the
quotient by all the previously chosen Chow families. The construction
terminates with rational-chain connectedness, and the classes of the
chosen families are numerically independent.

\begin{proposition}\label{prop:minimal-unsplit}
Under the assumption $\rho_X(\iota_X-1)\geq n-1$, every family in
this minimal sequence has anticanonical degree at most $\iota_X+1$
and is unsplit.
\end{proposition}

\begin{proof}
Lemmas~\ref{lem:short-degrees} and~\ref{lem:unsplit}, together with
Proposition~\ref{prop:chains}, give unsplit families
$V_1,\ldots,V_{\rho_X}$ such that
\[
 N_1(X)_{\RR}=\bigoplus_{j=1}^{\rho_X}\RR[V_j],
 \qquad -K_X\cdot V_j\leq\iota_X+1,
\]
and every point starts a chain containing a member of each family.
By Lemma~\ref{lem:rc-connected}, these families rationally chain
connect $X$. Proposition~\ref{prop:chains} also shows that $V_1$ is
covering. Minimality of $W_1$ therefore gives
\[
 -K_X\cdot W_1\leq -K_X\cdot V_1\leq\iota_X+1<2\iota_X.
\]
It follows from Lemma~\ref{lem:unsplit} that $W_1$ is unsplit.

At any subsequent nonterminal step, Lemma~\ref{lem:horizontal}
provides a horizontal dominating family among the original $V_j$.
The minimal choice $W_t$ consequently satisfies
\[
 -K_X\cdot W_t\leq -K_X\cdot V_j\leq\iota_X+1<2\iota_X,
\]
so it too is unsplit.We now argue inductively. Suppose that
$W_1,\ldots,W_{t-1}$ have been chosen according to the construction,
with arbitrary choices whenever the minimum degree is attained by
more than one family, and that these families are unsplit.
We take the rational-chain quotient
\[
 q_{t-1}\colon U_{t-1}\longrightarrow Z_{t-1}
\]
associated with the Chow families of all the previously chosen
families $W_1,\ldots,W_{t-1}$. If $Z_{t-1}$ is a point, the
construction terminates. Otherwise, $q_{t-1}$ is proper and
surjective, $U_{t-1}\subset X$ is dense open, and $Z_{t-1}$ is an
integral quasiprojective variety of positive dimension.

The original families $V_1,\ldots,V_{\rho_X}$ remain unsplit and
rationally chain connect $X$. Lemma~\ref{lem:horizontal}, applied
to $q_{t-1}$ and these families, therefore gives an index $j(t)$
such that $V_{j(t)}$ is horizontal dominating for $q_{t-1}$.
It is thus one of the families among which the minimum defining
$W_t$ is taken. Consequently, every possible choice of $W_t$
satisfies
\[
 -K_X\cdot W_t
 \leq -K_X\cdot V_{j(t)}
 \leq \iota_X+1
 <2\iota_X.
\]
Lemma~\ref{lem:unsplit} shows that $W_t$ is unsplit, completing
the induction step.
\end{proof}

Proposition~\ref{prop:minimal-unsplit} completes the passage from the
independent correspondences to the minimal families used by Fujita.

\Needspace{20\baselineskip}
\section{The classification}\label{sec:classification}

We first retain an auxiliary integer $s$ as a lower bound for the
Picard number. This is the formulation used in Fujita's inductive
step. It allows the base of the resulting contraction to have
Picard number greater than $s-1$, and it is suited to induction on
dimension, simultaneously for all auxiliary parameters. We now prove Fujita's advanced generalized Mukai conjecture in its formulation with an auxiliary lower bound for the Picard number.

\begin{theorem}\label{thm:auxiliary}
Suppose that $s\geq1$ is an integer and
\begin{equation}\label{eq:auxiliary}
 \rho_X\geq s,\qquad
 \iota_X\geq\frac{n+s-1}{s}.
\end{equation}
For $s=1$, we have $X\simeq\PP^n$ or $X\simeq Q^n$.
For $s\geq2$, $X$ is isomorphic to one of
\[
\begin{gathered}
 (\PP^{\iota_X-1})^s,\qquad
 Q^{\iota_X}\times(\PP^{\iota_X-1})^{s-1},\\
 B_{\iota_X}\times(\PP^{\iota_X-1})^{s-2},\qquad
 I_{\iota_X}\times(\PP^{\iota_X-1})^{s-2},\\
 \PP^{\iota_X}\times(\PP^{\iota_X-1})^{s-1}.
\end{gathered}
\]
Here the quadric factor may have dimension two.
\end{theorem}

Our first task is to obtain an extremal ray whose contraction has
fibers of dimension at most $\iota_X-1$.
Proposition~\ref{prop:minimal-unsplit} controls the successive
minimal sequence described above, and Fujita's contraction
criterion converts precisely this unsplitness property into the
required geometric conclusion. We recall the criterion with its
auxiliary parameter denoted by $s$.

\Needspace{10\baselineskip}
\begin{quote}
\noindent
\textbf{Fujita's contraction criterion
\cite[Corollary 6.4]{Fujita}.}
\emph{Suppose that $X$ is a smooth complex projective Fano variety
of dimension $n$, and that an integer $s\geq2$ satisfies}
\[
 \rho_X\geq s,\qquad
 \iota_X\geq\frac{n+s-1}{s}.
\]
\emph{Let $W_1,\ldots,W_k$ be the minimal sequence obtained by
choosing a dominating family of least anticanonical degree and then,
at each nonterminal step, a horizontal dominating family of least
anticanonical degree for the rational-chain quotient by the
previously chosen Chow families. If every $W_t$ is unsplit, then
there is an extremal ray $R\subset\NE(X)$ such that every fiber $F$
of its contraction $\varphi_R\colon X\to X_R$ satisfies}
\[
 \dim F\leq\iota_X-1.
\]
\end{quote}

This criterion supplies the contraction needed for induction.
To determine its total space, we also need the classification in
the dimension of its base. We make this second input precise.
For positive integers $d$ and $t$, we write $\AGM^d_t$ for the
entire assertion of Theorem~\ref{thm:auxiliary} in dimension $d$
with auxiliary parameter $t$. Thus $\AGM^d_t$ concerns every
smooth complex projective Fano variety $Z$ satisfying
\[
 \dim Z=d,\qquad
 \rho_Z\geq t,\qquad
 \iota_Z\geq\frac{d+t-1}{t},
\]
and gives the classification stated there, including the case
$t=1$.

In Fujita's inductive step, the small-fiber contraction is a smooth
$\PP^{i-1}$-fibration, where $i=\iota_X$. Its base has dimension
$n+1-i$, and its Picard number drops by one. The next result
therefore requires $\AGM^{\,n+1-i}_{s-1}$. The lower-dimensional
classification identifies the base with a rational model. Fujita's
result includes the passage from this fibration to a vector bundle
projectivization and the classification of the possible total
spaces. We use this conditional classification result in the
following form.

\Needspace{10\baselineskip}
\begin{quote}
\noindent
\textbf{Fujita's inductive classification
\cite[Proposition 5.3(b)]{Fujita}.}
\emph{Suppose that $X$ is a smooth complex projective Fano variety
of dimension $n$ and pseudoindex $i$, and that an integer $s\geq2$
satisfies}
\[
 \rho_X\geq s,\qquad
 i\geq\frac{n+s-1}{s}.
\]
\emph{Assume that there is an extremal ray $R\subset\NE(X)$ such
that every fiber of its contraction has dimension at most $i-1$.
If $\AGM^{\,n+1-i}_{s-1}$ holds, then $X$ is isomorphic to one of
the five varieties in Theorem~\ref{thm:auxiliary}, with auxiliary
parameter $s$ and $\iota_X=i$.}
\end{quote}

\begin{proof}[Proof of Theorem~\ref{thm:auxiliary}]
We prove the assertion by induction on $n$, simultaneously for
all auxiliary integers $s$. Thus, at the induction step in
dimension $n$, we assume $\AGM^d_t$ for every positive integer
$d<n$ and every integer $t\geq1$.

If $n=1$, then $X\simeq\PP^1$, and $\rho_X\geq s$ forces $s=1$.
For every $n$, the case $s=1$ follows from Dedieu--H\"oring's
numerical characterization of projective spaces and quadrics
\cite[Theorem 1.3]{DH}. We may therefore assume $n\geq2$ and
$s\geq2$.

We put $r=\rho_X$ and $i=\iota_X$. The hypothesis
\eqref{eq:auxiliary} implies $r\geq s\geq2$, $i\geq2$, and
\[
 r(i-1)\geq s(i-1)\geq n-1.
\]
Proposition~\ref{prop:minimal-unsplit} applies and shows that
every family in the minimal sequence is unsplit. This is the
sequence occurring in Fujita's contraction criterion. Since
\eqref{eq:auxiliary} also gives its numerical hypotheses, the
criterion produces an extremal ray $R\subset\NE(X)$ whose
contraction has every fiber of dimension at most $i-1$.

We next verify the lower-dimensional hypothesis in Fujita's
inductive classification. We set
\[
 n'=n+1-i.
\]
Theorem~\ref{thm:GM} gives $r(i-1)\leq n$, so
\[
 n'=n-(i-1)\geq(r-1)(i-1)\geq1.
\]
Moreover, $i\geq2$ implies $n'\leq n-1<n$.
Thus $n'$ is a strictly smaller positive dimension. Our induction
hypothesis applies with $d=n'$ and $t=s-1\geq1$, and proves the
entire statement $\AGM^{n'}_{s-1}$.

We also spell out why this is the classification needed for the
base in Fujita's inductive step. The base $Y$ supplied by that
step is a smooth projective Fano variety with
\[
 \dim Y=n',\qquad
 \rho_Y=r-1\geq s-1,\qquad
 \iota_Y\geq i.
\]
Its numerical conditions satisfy
\[
\begin{aligned}
 (s-1)(\iota_Y-1)
 &\geq(s-1)(i-1)\\
 &=s(i-1)-(i-1)\\
 &\geq n-i=n'-1.
\end{aligned}
\]
Hence $Y$ lies in the range of $\AGM^{n'}_{s-1}$. 

We have verified all the hypotheses of Fujita's inductive
classification: $X$ satisfies \eqref{eq:auxiliary}, the required
extremal ray exists, and $\AGM^{\,n+1-i}_{s-1}$ holds by induction.
Its conclusion gives exactly the five models in the statement.
\end{proof}

\begin{proof}[Proof of Theorem~\ref{thm:main}]
Lemma~\ref{lem:short-degrees} gives $\iota_X\geq2$ and
$\delta\in\{0,1\}$. If $\delta=0$, the equality case of
Theorem~\ref{thm:GM} gives
\[
 X\simeq(\PP^{\iota_X-1})^{\rho_X}.
\]

Suppose that $\delta=1$ and $\rho_X=1$. Then $n=\iota_X$.
The numerical characterization in \cite[Theorem 1.3]{DH} gives
$X\simeq\PP^n$ or $X\simeq Q^n$. The projective-space case has
pseudoindex $n+1$ and is excluded. Thus $X\simeq Q^n$.
Since $Q^2$ has Picard number two, we have $n\geq3$.

We finally suppose that $\delta=1$ and $\rho_X\geq2$.
We apply Theorem~\ref{thm:auxiliary} with $s=\rho_X$.
The model $(\PP^{\iota_X-1})^{\rho_X}$ has dimension
$\rho_X(\iota_X-1)$ and hence defect zero, so it does not occur.
If $\iota_X=2$, the quadric model is
\[
 Q^2\times(\PP^1)^{\rho_X-1}
 \simeq(\PP^1)^{\rho_X+1}.
\]
Its actual Picard number is $\rho_X+1$, so it too is excluded.
The remaining models are precisely those in part~\textup{(iii)}.
\end{proof}

\begin{corollary}\label{cor:exception}
Under \eqref{eq:auxiliary}, the only case with $\rho_X>s$ is
\[
 \iota_X=2,\qquad n=\rho_X=s+1,\qquad
 X\simeq(\PP^1)^{s+1}
       \simeq Q^2\times(\PP^1)^{s-1}.
\]
\end{corollary}

\begin{proof}
There is no such case in dimension one. For $n\geq2$, we have
$\iota_X\geq2$, and Theorem~\ref{thm:GM} gives
\[
 n-1\leq s(\iota_X-1)
 \leq\rho_X(\iota_X-1)\leq n.
\]
If $\rho_X>s$, then
\[
 1\leq(\rho_X-s)(\iota_X-1)\leq1.
\]
Both factors are positive integers, so $\rho_X=s+1$ and
$\iota_X=2$. The preceding inequalities then give $n=s+1$.
The equality case of Theorem~\ref{thm:GM} identifies
$X\simeq(\PP^1)^{s+1}$.
\end{proof}
\subsection*{AI disclosure} 
On 24 September 2026, OpenAI posted the proof of the generalised Mukai conjecture. We read this manuscript carefully by asking many questions to Chatgpt 6 and realised the techniques there can be used to prove Fujita's advanced generalised Mukai conjecture. We take full responsibility for the correctness and content of the paper.

\bigskip
\noindent\begin{minipage}{\textwidth}
\small
\textsc{Duo Li}\\
School of Mathematics, Sun-Yat Sen University,\\
Zhuhai, Guangdong 519082, PR China\\
\textit{Email:} \href{mailto:liduo5@mail.sysu.edu.cn}{\texttt{liduo5@mail.sysu.edu.cn}}
\end{minipage}

\end{document}